\documentclass[11pt]{amsart}
\usepackage{amssymb}
\voffset- .3cm
\numberwithin{equation}{section}

\newcommand{\Z}{\mathbb{Z}}

\usepackage{amsmath,amssymb}
\usepackage{picins,slashbox,
  graphicx,verbatim,longtable,
  epic,eepic,stmaryrd,rotating}
\usepackage[all]{xy}

\theoremstyle{plain}

\theoremstyle{definition}

\theoremstyle{remark}

\newlength{\standardunitlength}
\catcode`\@=11
\long\def\@makecaption#1#2{%
    \vskip 10pt
    \setbox\@tempboxa\hbox{
      \small{#1: }\ignorespaces #2}%
    \ifdim \wd\@tempboxa >\captionwidth {%
        \rightskip=\@captionmargin\leftskip=\@captionmargin
        \unhbox\@tempboxa\par}%
      \else
        \hbox to\hsize{\hfil\box\@tempboxa\hfil}%
    \fi}

\newdimen\@captionmargin\@captionmargin=2\parindent
\newdimen\captionwidth\captionwidth=\hsize
\catcode`\@=12

\newlength{\globalparindent}
\begin{document}

\vskip 8cm

\centerline {\bf The Finiteness Result for Khovanov Homology and Localization in Monoidal Categories.}

\vskip 1cm

\centerline {\bf Nadya Shirokova.}

\vskip 1cm

\centerline {\bf Abstract.}

\vskip .5cm 

  In [S1] we constructed the local system of Khovanov complexes on the Vassiliev's space of knots and extended it to the singular locus. In this paper we introduce the definition of the homology theory (local system or sheaf) of finite type and prove the first finiteness result: the Khovanov local system restricted to the subcategory of knots of the crossing number at most $n$ is the theory of type $\leq n$. This result can be further generalized to the categorification of Birman-Lin theorem [S2].

 \vskip 1cm

 \centerline {\bf Contents.}

\vskip 1cm 

1. Introduction.
  \vskip .3cm
2.  The Geometric Interpretation of the Cone of the Wall-Crossing Morphism.
  \vskip .3cm

3. Jones Polynomial and Invariants of Finite Type. 
  \vskip .3cm
4. The Grothendieck Group and the Reduced Khovanov Homology .
  \vskip .3cm
  
5. $\bigsqcup S^1$- stable Homotopy.  Category of Khovanov Spectra.
 \vskip .3cm 
6. The Poincare Polynomial.
   \vskip .3cm
7. On the Algebraic  Definition of Finiteness.
   \vskip .3cm
8. The Filtration.
   \vskip .3cm
9. The  Finiteness Result.
   \vskip .3cm

10. Further Directions.
\vskip .3cm
11. Bibliography.

 \newpage
 \centerline {\bf 1. Introduction.}

\vskip 1cm

 In [S2] we outlined a program of classification of homological knot theories, such as Khovanov's
 categorification of the Jones polynomial [Kh], Ozsvath-Szabo categorification of Alexander polynomial
 [OS] and Khovanov-Rozansky homology [KR].  This program can be further generalized to 3 and 4-manifolds.
 
\vskip .2cm
  
  Namely, we want to classify the following functors (embedded TQFT's):
  
    \[
\xymatrix{
\mathrm \mathcal Spaces \ar[rrr] ^-{functor}_{\mbox{\footnotesize(e.g. Khovanov)}} &&& \mathrm  \mathcal Triang.Cat}\\ \]

which behave in a prescribed way (via the wall-crossing morphisms) under cobordisms. By Spaces we understand the moduli spaces of manifolds, including the singular ones (e.g. Vassiliev's space of knots).

\vskip .2cm
        We  consider a knot homology theory as a local system, or a constructible
 sheaf on the space of all objects (knots, including singular ones), extend this local system to the singular locus and introduce the analogue of the "Vassiliev derivative" for categorifications.

  \vskip .3cm
  
   We get the correspondence between knots (possibly singular) and objects of the triangulated category (spectra), satisfying the exact triangle relations.
   \vskip .2cm
   
   Our classification  is given by the "type" of the theory, the homological condition on its extension to the strata of the discriminant, similar to the Vassiliev's classification of knot invariants.
   
   Recall that by Vassiliev, the knot invariant is of finite type $n$, if for any selfintersection of the discriminant of the space of knots of codimension $n+1$, the alternated sum of the invariants of knots from $2^{n+1}$ adjacent chambers is zero. 
  
  \vskip .3cm 
  The main example of this paper is the Khovanov homology, however  our results can be generalized
 to Khovanov-Rozansky theory.
    Recall that    M.Khovanov [Kh] categorified the Jones polynomial, i.e. he found a homology theory, the
  Euler characteristics of which equals the Jones polynomial.  From the diagram of the knot he
  constructs a bigraded complex, associated with this diagram, using 0 and 1- resolutions of the knot crossings.
 
\begin{center}
\includegraphics[scale=1]{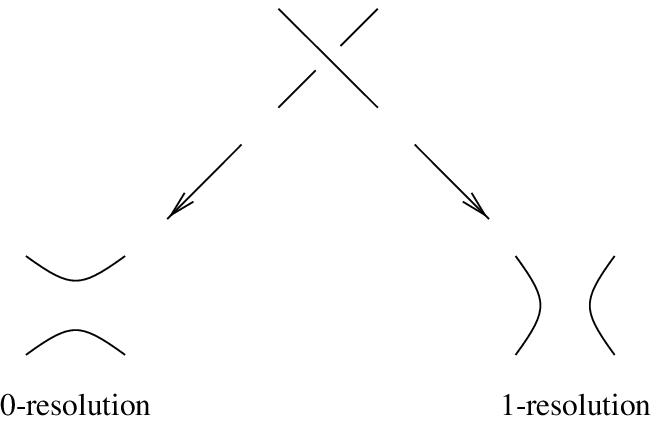}
\end{center}

  The  complex  becomes the sum of the tensor products of the vector space V, where the
  homological degree is given by the number of 1's in the complete resolution of the knot diagram.
  In [S1] we constructed the  local system of Khovanov complexes on the space of knots, the wall-crossing morphisms for the local system  and introduced the definition of the Khovanov homology of a singular knot.
 
  The discriminant of  the space of knots corresponds to knots with transversal self-intersection, it is an algebraic hypersurface which is cooriented, so that all cobordisms between knots are directed, i.e. moving between chambers we change overcrossing to undercrossing by passing through a knot
  with a single double point. We  extended the Khovanov local system to the discriminant
   by the cone of a wall-crossing morphism [S1]:

\vskip .5cm

{\bf Definition 1 [S1].} The Khovanov homology of the singular knot (with a single double point) is a bigraded complex 
$$
\renewcommand\arraystretch{1}
X^{\bullet}\oplus Y^{\bullet}[1] \quad\mathrm{with \ the \ matrix \ differential}\quad
d_{C_\omega}=\left(\!\!\!\begin{array}{cc}
d_X & \omega \\ 0 & d_Y[1] \end{array}\!\!\right),
$$
where $X^{\bullet}$ is Khovanov complex of the knot with overcrossing, $Y^{\bullet}$ is the Khovanov complex of the knot with undercrossing and $\omega$ is the wall-crossing morphism. 
\vskip .2cm
 
  In this paper we give the geometric interpretation of the homology of  the cone of the wall-crossing morphism:

\vskip .2cm
 {\bf Definition 2.} The Khovanov homology of the singular knot  $K$ of n crossings(with kth single double point ) is a bigraded complex 
$$
\renewcommand\arraystretch{1}
C^{\bullet}\oplus C^{\bullet}[2] \quad\mathrm{with \ the \ matrix \ differential}\quad
d_{C_\omega}=\left(\!\!\!\begin{array}{cc}
d_C & 0 \\ 0 & d_C[2] \end{array}\!\!\right),
$$
where $C^{\bullet}$  is the Khovanov complex of the knot of (n-1) crossings where  kth double point of $K$ is given 1-resolution.

\vskip .3cm
  This implies the geometric definition of the local system  being of finite type n:
      \vskip .5cm
  {\bf  Definition (G)}. The local system of Khovanov complexes, extended to the discriminant of the space of manifolds via the cone of morphism, is a {\bf local system of
  order n} if for any selfintersection of the discriminant of codimension $n$, its $n$'s cone is 
  quasiisomorphic to
  
  $$C_m^0 X^{\bullet}\oplus  C^{1}_{n}X^{\bullet}[2] ...\oplus C_n^n X^{\bullet}[2n]$$
  
  where $X^{\bullet}$ is the Khovanov complex of the disconnected sum of unknots.
  
\vskip .2cm

  In the triangulated category to every morphism between complexes there corresponds an object (up to isomorphism),  the mapping cone, which fits into an exact sequence.  By assigning cones of wall-crossing morphisms to singular knots we get the structure of the triangulated category for the Khovanov's sheaf.    We observe, that the constructed category has the monoidal symmetric  structure: the Khovanov complex of the disconnected sum of knots is the tensor product of the corresponding complexes.
  
   To give the definition of the local system of finite type and to prove the main result we construct  the category of the Khovanov spectra. 
    We  stabilize this triandulated monoidal category by taking disconnected sums of a knot with the collection of  circles, i.e. by taking tensor products of the sheaf with the complexes, corresponding to the disconnected sums of circles. (Tensoring with the Khovanov complex of a circle can be viewed as the suspension.) We use the reduced version of Khovanov homology [Kh1]. As the result  of the localization we get the category of Khovanov spectra $\tilde {\mathcal D}$ which  is again the symmetric monoidal category [Vo].
 \vskip .4cm

   Next we construct the sequence of derived categories, a filtration, obtained by factorization
   over the ideals, generated by the Khovanov complexes of the the special form, supported on the strata of the discriminant:
   
        $$ \tilde{\mathcal D}=\tilde{\mathcal D_{\infty}} \supset ...\tilde{\mathcal D_n} \supset  \tilde{\mathcal D_{n-1}} \supset ... \supset  \tilde{\mathcal D_1} $$ 
    
 \vskip .3cm
 
  This filtration is an analogue of the filtration in the theory of invariants of finite type. We show that in these subcategories the Khovanov theory satisfies the finiteness condition introduced in [S1]. 
  \vskip .5cm   
   The algebraic, Vassiliev-type definition of finiteness now becomes as follows.   Consider a new invariant,  an additive functor from complexes to the 2-torsion in their cohomology:
  
  $$\mathcal T_2: C \in Ob(\mathcal D) \rightarrow {Tor}_2( H^*(C))$$

  If the local system of Khovanov's complexes $CKh$ has  $\mathcal T^n_2=0$, i.e. $\mathcal T_2(H^*(CKh))|_{D_n}=0$ everywhere on $D_n$ - the codimension $n$ of the discriminant, but  is not quasiisomorphic to 
  $(Z^2,0,Z^2)$ (to eliminate the case of the Hopf link),  form a factor-category in a sense of Verdier $\mathcal D_n = \mathcal D/ \mathcal I_n$, where the category $\mathcal I_n$ is {\bf supported on the codimension n of the discriminant}. 
 
  \vskip .5cm

    Now can give an algebraic (torsion) definition, similar to the original one of Vassiliev (triviality  or acyclicity of complexes in codimension $n$) [S1].
        \vskip .3cm 

  {\bf Definition (T)}. The local system is of finite type n if for any codimension $n$ selfintersection of the discriminant its nth cone is not quasiisomorphic to 
  $(Z^2,0,Z^2)$, and has  torsion-free homology  (i.e. the image of  $\mathcal T_2$
  is zero).
    \vskip .3cm

  The Vassiliev-type definition becomes:
 
   \vskip .3cm

 {\bf Definition (V)}. The local system of Khovanov complexes is of finite type $n$ if for any codimension $n$ selfintersection of the discriminant the nth cone is zero in $\tilde{\mathcal D_n}$ but not in
 $\tilde{\mathcal D_{n+1}}$.

  \vskip .3cm

 The main results of this paper  is as follows:
   \vskip .3cm

{\bf Theorem  1}. Restricted to the subcategory of knots with the crossing number  at most $n$, $n\geq 3$, Khovanov local system is of finite type $\leq n$.

\vskip.3cm

 We will further  generalize this result [S2] and get the "categorification of Birman-Lin theorem" [S2].

  \vskip .5cm
  {\bf Acknowledgements}.  I want to thank P. Deligne, G. Carlsson,  O. Viro and A. Voronov   for  useful discussions  and  Stanford University for their hospitality.  

\newpage

 \centerline {\bf 2. The geometric interpretation of the cone of the wall-crossing morphism. }

\vskip 1cm
  In this paragraph we give the geometric interpretation of the cone of the wall-crossing morphism
  for Khovanov homology, by raising the skein relation for the Jones polynomial to the level of complexes.

   First we recall the definition of the wall-crossing morphism $\omega_k$, correponding to the crossing
   change for the kth point in the knot projection. 
    
  \vskip 1cm
 
\begin{equation*}
\xymatrix@C+0.5cm{\omega_k: A_0^\bullet (k)
\ar[r]^-{ Id} & B_0^\bullet[1] (k)  \\
\omega_k:A_1^\bullet(k)[1] \ar[r]^{ \emptyset} & B_1^\bullet[1](k) }
\end{equation*}
 \vskip 1cm
 
 Where $A^\bullet$ is the Khovanov complex, corresponding to the knot with the kth overcrossing and
 $B^\bullet[1]$ to the knot with kth undercrossing.
 
  When we change kth overcrossing to undercrossing, 0 and 1-resolutions are exchanged , so $A^\bullet = A^\bullet_0(k) \oplus A^\bullet_1(k)$, $B^\bullet[1]=B^\bullet_0(k)\oplus B^\bullet_1[1](k)$, and for every k we get morphism  $\omega_k$.
  
    \vskip .5cm
    According to the Definition 2 to show that the local system is of finite type (n-1), we have to prove that
   for any point of selfintersection of the discriminant of codimension n the corresponding n-cones are
   acyclic complexes, or that the convolution of n dimensional hypercube is acyclic.

   \vskip.3cm
   
   We first observe the following simple properties of the Khovanov complex:
    \vskip.2cm 
    
    {\bf Proposition 1}. For any commutative n-dimensional cube , with $2^n$ complexes of length $(n+1)$  at its vertices, the last complex (via coorientation) is
    the the dual of the first one.
     \vskip.2cm 
    {\bf Proof}. Notice, that by changing all overcrossings to undercrossings on the diagram projection we exchange 0-resolutions of the diagram  to 1-resolutions and move from the knot to its mirror. 
 One can easily see that the Khovanov complex will be dualized.
   \vskip .4cm 
      Now we want to study the restrictions of our local system.
       Recall that  a $\bf subcategory$ of a category C is a category S whose objects are objects in C and whose arrows  $ f:A\rightarrow B$ are arrows in C (with the same source and target). 
   
    For the category of knots $\mathcal K$  we can define a sequence of subcategories $\mathcal K_n$ . 
    Objects
    of $\mathcal K_n$ are knots with a crossing number at most n. 
   \vskip.2cm

      We refer to [GM] for the definition and properties of Postnikov towers and convolutions.
       Here we show that for the restriction of the Khovanov local system to the subcategory 
       $\mathcal K_n$ of knots with at most n crossings, complexes forming equators of the hypercubes
       will fit into a sequence.      
     \vskip.5cm 
      
     {\bf Proposition 2}.  If n complexes $X^\bullet,....,Z^\bullet$ form the equator of the n-hypercube 
     in the restriction of the Khovanov local system to $\mathcal K_n$, the subcategory of knots with at most n crossings, 
     then $X^\bullet,....,Z^\bullet$ and the wall-crossing morphisms $\omega$
     fit into a sequence:
         \vskip .5cm

    
      $$\rightarrow Z^\bullet [-n+1] \rightarrow^u X^\bullet \rightarrow ^{\omega}...\rightarrow ^
      {\omega}Z^\bullet \rightarrow^u X^\bullet[n-1]$$
      
        \vskip .5cm 
      where $u$  is an isomorphism between $X^n$ and $Z^0$.
     \vskip.3cm 
     
     {\bf Proof}. 
      If $(n+1)$ complexes  $X^{\bullet},....,Z^{\bullet}$ form the equator in the n-dimensional 
      commutative hypercube,
      there is a connecting map $u$ from the last complex into the first one, shifted by [n-1], [GM]. Given
       Proposition 1 , the connecting map u is the isomorphism between $Z_0$ and $X_n$.
      
        \vskip 1cm

 {\bf Proposition 3}. The cone of the wall-crossing morphism between Khovanov complexes corresponding to $K$ and  $K'$, where the ith double point of the projection of $K$ is overcrossing and 
 the ith double point of $K'$ is undercrossing, is quasiisomorphic to
 $$X^{\bullet}\oplus X^{\bullet}[2]$$
 
  where $X^{\bullet}$  is the Khovanov complex of the knot of (n-1) crossings where  ith double point of the projection of $K$ is given 1-resolution.
   \vskip .5cm
 {\bf Proof}.  It is clear from the definition of $\omega_k$ and the way we defined the local system, 
  that the wall-crossing morphism corresponding to the ith crossing will map isomorphically the parts of the complex, which have 0-resolution of the crossing (and will be quasiisomorphic to zero in the cone of the wall-crossing morphism).  The parts, that are mapped by zero, are the 1-resolutions of the k th crossing. In particular, we see that the last component of the complex
 contributes nontrivially to the homology  of the cone.
 
 \vskip .3cm
         
        So the homology of the cone of the wall crossing morphism is isomorphic to the homology of
        the knot with the ith intersection point given 1-resolution.
        
        This allows us to give the geometric version of the Definition 1:
        
        Let $D$ be the projection of the knot of n crossings with k th crossing being a doublepoint. We call
        $D$ the projection of the singular knot. 
    \vskip .5cm
 {\bf Definition 2.} The Khovanov homology of the singular knot  $K$ of n crossings(with kth single double point ) is a bigraded complex 
$$
\renewcommand\arraystretch{1}
C^{\bullet}\oplus C^{\bullet}[2] \quad\mathrm{with \ the \ matrix \ differential}\quad
d_{C_\omega}=\left(\!\!\!\begin{array}{cc}
d_C & 0 \\ 0 & d_C[2] \end{array}\!\!\right),
$$
where $C^{\bullet}$  is the Khovanov complex of the knot of (n-1) crossings where  kth double point of $K$ is given 1-resolution.

  \vskip 1.5cm

        On the level of the Euler characteristics one has:
        
        $$\chi (A^\bullet[1])= - \chi (A^\bullet[1])$$
        
        $$   \chi (C_f) = \chi (A^\bullet) -  \chi (B^\bullet)$$
        
        Since the Euler characteristics of Khovanov complex is the Jones polynomial, recall its skein
        ralation:

        $$q^{-1} J_{L_+} - qJ_{L_-} = (q^{1/2} - q^{-1/2}) J_{L_0}$$

        where $L_+$ is the knot with overcrossing , $L_-$ knot with undercrossing and $L_0$ is
        the knot, where the crossing point  is given 1-resolution.
    \vskip .5cm

 {\bf Proposition 4}. Let  $K$ be the knot with n crossings, denote it's projection $D_{i_1,...,i_n}$, where
 each index $i_n$ can have values $+$ for overcrossing, $-$ for undercrossing and $0$ for
 1-resolution of the crossing point. Then the m'th cone of the local system, for which $K$ is the "first"
 knot is given by formula (for simplicity we assume that the wall-crossings happen for the first m 
 indices):
  
 $$C_m^0 CKh(D_{0...0,i_{m+1},..i_n}) \oplus C^{1}_{m} CKh(D_{0...0,i_{m+1},..i_n})[2] \oplus...$$
 $$...\oplus C^{l}_{m} CKh(D_{0...0,i_{m+1},..i_n})[2l]\oplus ... \oplus CKh(D_{0...0,i_{m+1},..i_n})[2m]$$
 \vskip .5cm
 {\bf Proof}.  We will prove this formula by induction. The first cone was described in Proposition 3. In
 our notation it is given by formula $CKh(D_{0,i_2,..i_n}) \oplus CKh(D_{0,i_2,..i_n})[2]$. If we take the second cone, we get
 $$CKh(D_{0,0,i_3,..i_n}) \oplus CKh(D_{0,0,i_3,..i_n})[2] \oplus CKh(D_{0,0,i_3,..i_n})[2] \oplus CKh(D_{0,0,i_3,..i_n})[4]$$
 
 Suppose we proved the formula for m-1, let's show it is true for m.
 
\vskip .3cm

 The $m-1$st cone is found to be
 
  $$C_{m-1}^0 CKh(D_{0...0,i_{m},..i_n}) \oplus C^{1}_{m-1} CKh(D_{0...0,i_{m},..i_n})[2] \oplus...$$
 $$...\oplus C^{l}_{m-1} CKh(D_{0...0,i_{m},..i_n})[2l]\oplus ... \oplus CKh(D_{0...0,i_{m},..i_n})[2m-2]$$
 \vskip .2cm
 When taking the last mth cone, we will take pairwise cones of corresponding summands and can use Proposition 3 to show that we will be getting the components of the complex, described in Proposition 4.

  We just have to show that the coeffitients of the formula will be given by binomial coeffitients of proposition 4. This can be proved by using the identities for the binomial coeffitients.
For any n, m we use the formula: 
  \vskip .5cm
 
\begin{equation*}
\xymatrix@C+0.5cm{C^m_n = C^{m-1}_{n-1} \oplus C^m_{n-1}}
\end{equation*}
 \vskip .5cm

The geometric version of the definition of finiteness follows from the Definition 2:

    \vskip .5cm
  {\bf  Definition (G)}. The local system of Khovanov complexes, extended to the discriminant of the space of manifolds via the cone of morphism, is a {\bf local system of
  order n} if for any selfintersection of the discriminant of codimension $n$, its $n$'s cone is 
  quasiisomorphic to
  
  $$C_m^0 X^{\bullet}\oplus  C^{1}_{n}X^{\bullet}[2] ...\oplus C_n^n X^{\bullet}[2n]$$
  
  where $X^{\bullet}$ is the Khovanov complex of the disconnected sum of circles.

  In paragraphs 5,7 we will construct the sequence of derived categories, in which Definitions 2 and 4
  will become equivalent.

\vskip 1.5cm
 \centerline {\bf 3. Jones polynomial and invariants of finite type. }

\vskip 1cm	

The notion of the invariant of finite type was introduced by V.Vassiliev in 1989 as a filtration in the spectral sequence. Later	it
was  interpreted by Birman and Lin as a "Vassiliev derivative" and led to
the following skein relation.

If $\lambda$ be an arbitrary invariant of oriented knots in oriented space
with values in some abelian group $A$. Extend $\lambda$ to be an invariant of
$1$-singular knots $L_1$ (knots that may have a single singularity that locally
looks like a double point), using the formula

$$  \lambda(L_1)=\lambda(L_+)-\lambda(L_-)$$

where as before $L_+$ is the knot with overcrossing , $L_-$ knot with undercrossing.

Further extend $\lambda$ to the set  of $n$-singular knots $L_n$ (knots
with $n$ double points) by repeatedly using the skein relation.
    \vskip .2cm
{\bf Definition} We say that $\lambda$ is of type $n$ if its extension
to $(n+1)$-singular knots vanishes
identically. We say that $\lambda$ is of finite type if it is of type $n$ for
some $n$.
    \vskip .2cm	
    
    Let $\mathcal L_n$ be invariants of knots (with values in Q) of order $\leq n$, then
    $\mathcal L_n / \mathcal L_{n-1}$	 invariants of knots of order exactly n.

  Let $\mathcal L$ - formal linear combinations of knots , then the singular knot is a linear combination
  of $2^n$ terms.
  
  Let $\mathcal L^n$ - subspace of $\mathcal L$, generated by knots with n double points.
     \vskip .3cm 
  {\bf Fact 1}.   $\mathcal L^{n+1} \subset \mathcal L^n$ . 		
      \vskip .3cm
  
   {\bf Fact 2.} Spaces $\mathcal L_n / \mathcal L_{n-1}$	 and  $\mathcal L^n / \mathcal L^{n+1}$
   are dual to each other.	 
      \vskip .3cm		
Birman and Lin (1993) showed that substituting the power series for $e^x$ as the variable in the Jones polynomial yields a power series coefficients of which  are Vassiliev invariants:
     \vskip .3cm	
{\bf Theorem [BL]} Let  $K$ be a knot and $J_t(K)$ be its Jones polynomial. let $U_k(x)$ be obtained from $J_t(K)$  by replacing the variable $t$ with $e^x$. Express $U_k(x)$ as power series in$x$:

$$U_k(x) = \sum u_i(K)x^j$$

then $u_0(K) = 1$ and each $(K)\geq 1$ is a Vassiliev invariant of order i.
 \vskip .3cm

 Their result implies in particular that  the values of the Jones polynomial are not of finite type, but the values of the truncations are.
 \vskip .2cm
 
  Note, that there is another normalization of the Jones polynomial, when it is considered not as a
  polynomial in $(q + q^{-1})$, but in $(q-1)$. This renormalization implies the Jones polynomial, determined by the skein relation
  
  $$q^2 J_{L_+} - q^{-2} J_{L_-} = (q - q^{-1}) J_{L_0}$$
 \vskip .3cm 
  and normalized by $J(\bigcirc)=1$.
   \vskip .3cm	
   The Birman-Lin theorem holds for this renormalized polynomial in $q$ without substituting $e^x$.
  
  The Khovanov theory, categorifying this polynomial, is called the reduced Khovanov homology.

 Since  the Euler characteristics of the cone of the morphism is just the difference of the values of the
 Jones polynomials of knots from adjacent chambers,, our local system is not expected to be of finite type, however, as we
 will see later, it's restrictions can be of finite type.

\vskip 1.5cm

 \centerline {\bf 4. The Grothendieck group and the reduced Khovanov homology .}

\vskip 1cm

 In this paragraph we  start to  establish the correspondence between Vassiliev-type definition of finiteness given in [S1]  and the geometric one of section 2.
 . We introduce the Grothendieck group of Khovanov homology to be able to factorize  over subgroups , generated by the
 Khovanov  homology of the collection of circles and their shifts.

To construct the Grothendieck group of a commutative monoid M, one forms the Cartesian product

   $$ M \times M$$

The two coordinates  represent first and second part:

   $$ (a, b)$$

which corresponds to

    $$a - b$$

Addition is defined as follows:

  $$  (a, b) + (c, d) = (a + c, b + d)$$
\vskip .5cm
Next we define an equivalence relation on $M\times M$:  (a, b) is equivalent to (c, d) if, for some element k of M if $ a + d + k = b + c + k$. It is easy to check that the addition operation is compatible with the equivalence relation. The identity element is now any element of the form (a, a), and the inverse of (a, b) is (b, a).

The Grothendieck group can also be constructed using generators and relations: denoting by (Z(M),+) the free abelian group generated by the set M, the Grothendieck group is the quotient of Z(M) by the subgroup generated by $ \{a+ b - (a+b)\mid a,b\in M\}$.

The Grothendieck group
of an abelian category $\mathcal M$ is an abelian 
group with generators $[M],$ for all objects and
relations $[M_2]=[M_1]+[M_3]$ for all  exact sequences: 

$$ 0 \rightarrow M_1 \rightarrow M_2 \rightarrow M_3 \rightarrow 0$$

The Grothendieck group of the category $\mathcal K(\mathcal M)$ of bounded
complexes  up to chain homotopies is an abelian
group with generators $[M],$ and 
 relations $[M[1]]= - [M]$ 
and $[M_2]= [M_1]+[M_3]$ for 
all  exact sequences of complexes as above for all components of the complexes.

The Grothendieck group of a triangulated category 
$\mathcal T$ is an abelian group with generators $[M],$ for all objects
$M$ of $\mathcal T$ and relations $[M[1]]=-[M]$ and $[M_2]=[M_1]+[M_3]$ 
for all distinguished triangles 

$$ ...\rightarrow M_1 \rightarrow M_2 \rightarrow M_3 \rightarrow M_1[1] \rightarrow ...$$

It is easy to see that the 
Grothendieck group of the bounded derived category $D^b(\mathcal M)$ 
is isomorphic to the Grothendieck groups of $\mathcal K(\mathcal M)$ and 
$\mathcal M$.

The Grothendieck group was originally introduced for the study of Euler characteristics.

  Now we want to apply the above constructions to the derived category of Khovanov complexes.
  
   As it was shown in the original Khovanov paper, if $D_1,D_2$  are diagrams of oriented links $L_1,L_2$ then for 
a diagram $D_1\sqcup D_2$ of the disjoint union $L_1\sqcup L_2.$ 
there is an isomorphism of cochain complexes 
$$
C(D_1\sqcup D_2) = C(D_1) \otimes  C(D_2)$$
of free graded abelian groups. 
\vskip 1cm

From the K\"unneth formula Khovanov derives the following formulas for the cohomology of the disjoint union:
\vskip .5cm
 
{\bf Proposition.} There is a short split exact sequence of cohomology 
groups 

\vskip .3cm
$$0  \to   \oplus_{i,j\in Z}(  H^{i,j}(D_1)\otimes H^{k-i,m-j}(D_2))
\to H^{k,m} (D_1\sqcup D_2) \to 
 \oplus_{i,j\in \Z} 
 Tor_1^{\Z}( H^{i,j}(D_1), H^{k-i+1,m-j}(D_2))\to 0$$
\vskip .5cm

{\bf Corollary.} 
 For each $k,m \in Z$ there 
  is an equality of isomorphism classes of abelian groups 
  
  \vskip .3cm
 $$  H^{k,m} (L_1\sqcup L_2) = 
   \oplus_{i,j\in \Z}( H^{i,j}(L_1)\otimes  H^{k-i,m-j}(L_2)) 
  \oplus_{i,j\in Z} 
 Tor_1^{\Z}(H^{i,j}(L_1), H^{k-i+1,m-j}(L_2)) $$
\vskip .5cm

Over $Q$ these formulas will imply that the derived category of Khovanov complexes form a 
 tensor category which has a monoidal structure.

 Given the disconnected sum of two knots $K_1, K_2$ one gets a tensor product of Khovanov groups:

 $$H^{a,b} (K_1) \otimes H^{c,d} (K_2)$$

 We form the Grothendieck group as follows:
 
 Khovanov homology are the bigraded groups, so we will be taking sums over all products
 with fixed bigradings $a+c, b+d$:
 
 $$H^{a,b} \otimes H^{c,d} \rightarrow H^{a+c,b+d}$$
 
 In the original Khovanov's paper [Kh] the homology of a loop is 
 $$ H( \bigcirc) = Z\{-1\} \oplus Z\{1\}$$
 which corresponds to the normalization of the Jones polynomial, s.t.
  $$J( \bigcirc)= q+q^{-1}$$
 
 Then the Khovanov homology of the disconnected sum of m circles is 
 $$H (\bigcirc ^ m) = ({Z\{-1\} \oplus Z\{1\}})^{\otimes m}$$
  Thus in the Grothendieck group the identity will be formed by the elements
 of the form $H^{0,j}$ and after factorizing by the disconnected sum of circles and their shifts we get the new version of Khovanov homology $'Kh$  we will  get the identity:
 
 $$'H^{i,j}(D) = H ^{i+j,j}(D)$$

  The version of the Khovanov homology which was introduced in [Kh1] takes care of the above problem.
  It categorifies the Jones polynomial, determined by the skein relation
  
  $$q^2 J_{L_+} - q^{-2} J_{L_-} = (q - q^{-1}) J_{L_0}$$
 \vskip .3cm 
  and normalized by $J(\bigcirc)=1$.
  \vskip .5cm
  
  It is called the reduced homology and is defined as follows:
    \vskip .5cm
  
   Let $ \mathcal A = Q[X]/(X^2)$ is the base ring and $H^{i,j}(D)$ is the complex of finite-dimensional 
   $Q$-vector spaces. Khovanov constructs a map of complexes  $ \mathcal A \otimes C(D) \rightarrow C(D)$ via  a geometric construction: choose a segment of the knot diagram $D$ that doesn't contain crossing, place an unknotted circle next to it and consider cobordism, which merges circle into $D$. Ridemeister move, which happens away from this cobordism, induces a chain homotopy equivalence
   between complexes of  $ \mathcal A$-modules. It also establishes a bijection between $(1,1)$-tangles
   and oriented links with marked component.
 \vskip .5cm  

   If $\mathcal Q  = \mathcal A/X  \mathcal A $ is onedimensional representation of  $ \mathcal A$, then 
   the reduced complex is defined as
   
   $$ \tilde {C(D)} = C(D) \otimes_{ \mathcal A }  \mathcal Q $$
   
   and its homology is the reduced homology of $D$. The analogous construction can be carried out over integers.
   
   \vskip .5cm
   
   In our case we mark any arc of the "first" knot, which will become a circle after making 1-resolutions of
   all crossing point. By the cobordism construction all other circles of the resolved link can be merged to 
   this component.

 \vskip .3cm
 
  After we made a derived category of Khovanov complexes into a "group", the first guess of how to match the geometric and the Vassiliev-type definitions would be: factorize the category by the complexes corresponding to the disconnected sums of circles. However, disconnected sums with unknots will give acyclic complexes  and the category will become trivial.

  {\bf Note}. Quotients in the DG categories were studied by V. Drinfeld [D]. However, these are not
  the quotients that we would like to consider, since we don't want homology of the unknot to be zero.
  \vskip .3cm 
  We will do the factorization in two steps: first we construct a stable category, in which the disconnected sum with a circle is viewed as a suspension, so that the homology of the disconnected sum of a knot and a circle would be quasiisomorphic to the homology of the knot. (In that case we won't have to remmember how many circles our theory decomposed into).
  
  Next we form the filtration, on the factors of which the geometric definition will become equivalent to the Vassiliev's one.
  
 We will be using  the reduced version of Khovanov homology.

  \vskip 1cm

 \centerline {\bf 5. $\bigsqcup S^1$- stable Homotopy.  Category of Khovanov Spectra.}
 
  \vskip 1cm
   We introduce a new derived category $\tilde {\mathcal D}$ in which the reduced Khovanov homology of the disconnected sums of circles, discussed in the previous paragraph, is quasiisomorphic to the Khovanov homology of a circle.
   
     This can be understood  by considering a linear equivalence relation , generated by
  an object - unknot $u$: 
  $$X \sim Y  \Leftrightarrow X \otimes u_1 = Y \otimes u_2$$
  where $u_1$ and  $u_1$ are disconnected sums of unknots. 
   \vskip .3cm	
  What we will construct in this section can be viewed
  as a version of {\bf Spanier-Whitehead category}, objects of which are called ${\bf spectra}$.

 This construction can be carried out in any
 {\bf symmetric monoidal} category and it was proved  that the result of this localization is again a symmetric monoidal category, e.g. [Vo]. 
  
  
    Recall the definition of the monoidal category:
  
\vskip .5cm

 {\bf Definition}. A  {\bf monoidal category} (or tensor category) is a category $\mathcal M$ equipped with
\vskip .2cm
    1) a binary functor $\otimes \colon \mathcal M\times\mathcal M\to\mathcal M$  called the tensor product,
    
    2) an object I called the unit object,
    
    3) three natural isomorphisms subject to certain coherence conditions expressing the fact that the 
    tensor operation  $\otimes$ is associative, i.e. there is a natural isomorphism , called associativity, with components  $ \alpha_{A,B,C} \colon (A\otimes B)\otimes C \to A\otimes(B\otimes C)$,
          $\otimes$ has I as left and right identity: there are two natural isomorphisms ,  with components
 $\lambda_A \colon I\otimes A\to A $ and $\rho_A \colon A\otimes I\to A.$
 
 The coherence conditions for these natural transformations, given by penthagon diagrams which commute for all objects $A,B,C,D \in \mathcal M$ (see Fig. 1 on the next page).
\vskip 1cm
The tensor category, which we are considering is also {\bf braided}, i.e. it is equipped with the
braiding isomorphism $$ \gamma_{A,B} : A \otimes B \rightarrow B \otimes A$$

A {\bf symmetric} monoidal category is a braided monoidal category whose braiding satisfies
 $\gamma_{B,A} \gamma_{A,B}  = 1_{A \otimes B}$.

  \vskip .3cm

  It is obvious that the derived category of Khovanov complexes satisfies all these conditions and is a symmetric monoidal category [Kh2]. 

       \vskip .2cm
     {\bf Definition.} Suppose $\tilde {\mathcal D}$ satisfies the Definition 1 for {\bf countable} coproducts.
     Let

\begin{equation*}
\xymatrix@C+0.5cm
{X_0\ar[r]^-{ j_1} & X_1 \ar[r]^-{ j_2} & X_2 ...}
 \end{equation*}

 be a sequence of objects and morphisms in $\tilde{ \mathcal D}$. The homotopy colimit of the sequence
 denoted ${\mathcal Hocolim}_{i \to +\infty} (X_i)$ is by definition given, up to non-canonical isomorphism by the triangle:
 
        \vskip .5cm
 
\begin{equation*}
\xymatrix@C+0.5cm
 {\coprod X_i\ar[r]^-{ [1]-shift} &{ \coprod X_i}\ar[r] & {\mathcal Hocolim (X_i)} \ar[r] &\Sigma \{ \coprod X_i \} }
 \end{equation*}
 \vskip .5cm     
     
 here the shift map $[1]-shift$ is the infinite matrix:
  \vskip 1cm

\[
\left(
\begin{array}{ccccc}
1_{X_0}&j_1&0&1&...\\
0&1_{X_1}&j_2&0&...\\
0&0&1_{X_2}&j_3&...\\
0&0&0&1_{X_3}&...\\
...&...&...&...&...\\

\end{array}
\right)
\]

  \vskip 1cm 
\newpage
The coherence conditions:

\vfill
\begin{align*}
\xymatrix
{
 & (A\otimes B) \otimes (C\otimes D)
 \ar[ddr]^{\alpha_{A, B,C\otimes D}} & \\ \\
 ((A\otimes B)\otimes C)\otimes D
 \ar[uur]^{\alpha_{A\otimes B,C, D}}
 \ar[dd]^{\alpha_{A, B,C} \otimes 1_D}
 & &
 A\otimes(B\otimes(C\otimes D))\\ \\
(A\otimes (B\otimes C)) \otimes D
 \ar[rr]^{\alpha_{A, B \otimes C, D}}
 & &
A\otimes ((B\otimes C)\otimes D)
 \ar[uu]^{ 1_A \otimes \alpha_{ B , C, D}}
}
\end{align*}

  \vfill
  
    One can think of taking a disconnected sum with the circle as a suspension on Khovanov complex.
    Indeed, connect the circle to the knot by two strands, forming the overcrossing point in projection:
    \vfill
\begin{center}
\includegraphics[scale=1]{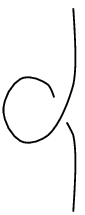}
\end{center}
\vfill
then the resolutions are
\vfill
\begin{center}
\includegraphics[scale=1]{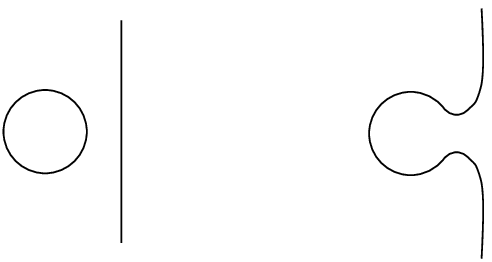}
\end{center}
\vfill
   
\newpage
    The resulting complex will be of length $(n+1)$ if the original was of length $n$ and will have
    shifted grading, like in the case of usual suspensions. Desuspending will correspond to eliminating
    fake loops as above.

    \vskip .2cm

     Then the Freudenthal theorem implies that
    
    $$[X,Y] \rightarrow [\Sigma X, \Sigma Y] \rightarrow [\Sigma^2 X, \Sigma^2 Y] \rightarrow...$$

    eventually stabilizes. The stable homotopy classes of maps from $X$ to $Y$ is above colimit.

   The resulting category is the Spanier-Whitehead category, or the  category of Khovanov spectra.
   
   From the category viewpoint the above construction implies, that we added another axiom to the
   the standard axioms of triangulated category:
     \vskip .2cm
     {\bf Definition.} Let $\alpha$ be an infinite cardinal. The triangulated category $\tilde {\mathcal D}$ is said to satisfy
     {\bf homotopy colimits} axiom if the following holds:
     
     For any set $\mathcal N$ of cardinality less than $\alpha$ and any collection 
     $\{ X_n, n\in \mathcal N \}$ of objects, the coproduct $\coprod X_n$ exists in $\tilde{ \mathcal D}$.

  Now we  define the new category $\tilde {\mathcal D}$, by modifying morphisms:
  
  \vskip .5cm
  
  {\bf Definition.} The category $\tilde {\mathcal D}$ of Khovanov {\bf  spectra} is as follows
  
 1) Objects of $\tilde {\mathcal D}$ = Objects of $ \mathcal D$

  2) Morphisms of $\tilde {\mathcal D}$:
  
 $$ \tilde Hom (X \rightarrow Y) = colim (Hom (X \rightarrow (Y + n (Unknots))))$$

 \vskip .5cm
 One can  define composition of morphisms.
  If $\mathcal H_n = Hom (X \rightarrow (Y + n (Unknots)))$, then one has to take the limit
  
  $$\mathcal H_1 \rightarrow \mathcal H_2 \rightarrow... \mathcal H_i \rightarrow...$$

  where $i$ denotes the number of unknots .

    \vskip .2cm 
    The fact, that the resulting category is triangulated and has nice properties  was proved by several
    authors [A], [Vo], in algebraic setting for example:
  \vskip .2cm    
  
 {\bf Theorem [Vo] }. Let $ (\mathcal D, \otimes , 1)$ be a symmetric monoidal category and $U$
 be an object, s.t. the cyclic permutation on $U \otimes U \otimes U$ equals identity in
  $\mathcal D[U^{-1}]$. Then there exists a symmetric monoidal structure $\otimes$ on $\mathcal D[U^{-1}]$.
   \vskip .5cm
 The above theorem  implies that in the category $\tilde {\mathcal D} = \mathcal D[U^{-1}]$ the disconnected sums
 of unknots are considered  an identity and the resulting category is again symmetric monoidal.
\vskip 1.5cm

\newpage
 \centerline {\bf 6. The Poincare Polynomial.}
 
  \vskip 1cm

  As we observed in paragraph 2, Proposition 4, the nth cone of the Khovanov local system, restricted to the subcategory of knots of crossing number n is a complex of length $2n+1$, with generators  (unknots)
 in the odd homological dimensions.
  \vskip .3cm
 {\bf Definition}. The {\bf Poincare polynomial} of the complex $X^{\bullet} =(X^0,X^1,....X^n)$ is the
 polynomial $P(X,t)=\Sigma{ \beta}_k t^k$, where ${\beta}_k = dim(X^k)$ is the kth Betti number.
 \vskip .3cm
We consider the ideals $ I_n$ generated by the complexes of the form:

 \[ C^\bullet (n) = \left (   \begin{array}{c}
   {C_n^n} Z \\ 0 \\ \\ 0 \\ {C^{i}_{n}} Z  \\ 0 \\ \\ 0 \\{C_n^0} Z  \end{array} \right) \]

 \vskip .5cm

In the next paragraph we will factorize the category of complexes by these complexes. Notice that
$C^\bullet (n)=(C^\bullet (1))^{\otimes n}$.
  \vskip .2cm
{\bf Lemma}. The Poincare polynomial of the complex $C^\bullet (n)$ equals
 ${(1 + t^2)}^n$ and its Euler characteristics is $2^n$.
   \vskip .2cm
 {\bf Proof}.     Obvious:
 $${(1 + t^2)}^n = \Sigma C^{i}_{n} t^{2i}$$ 
 
 $$\Sigma {C^{i}_{n}} = 2^n.$$
 
  \vskip .2cm
{\bf Note}.  After grading change, we can assume that the Poincare polynomial is equal to 
${(t^{-1} + t)}^n$.

For each n by the polynomial ${(t + t^{-1})}^n$ we  generate an ideal $ I_n$  in the ring of Laurent polynomials $Z[t, t^{-1}]$:

$$ I_n = Z[t, t^{-1}]/ {(t^{-1} + t)}^n$$
 \vskip .3cm
In the next chapter we are taking this factorization to the level of the derived category.
\vskip 1.5cm
\newpage

 \centerline {\bf 7.  On the Algebraic Definition of Finiteness.}
 
  \vskip 1cm
  
  The geometric definition of finiteness, which we gave in section 2 is sufficient to prove the main result of this paper, however, we would like to give an algebraic definition, similar to the original one of Vassiliev (that complexes become zero or acyclic after the extension to the  codimension $n$ of the discriminant). 
  
  However, we would like to understand if there are computational ways to determine that the local system decomposed into a collection of circles on the strata of the discriminant, i.e. how the TQFT "knows" about it?

 The idea is to factor out the free part of the homology.  
    \vskip .3cm    
 
   Recall the Theorem of A.Shumakovich [Su] and Asaeda-Przytycki [AP] regarding the torsion in Khovanov homology:
 \vskip .3cm   
   {\bf Theorem [Su]} The only alternating links that  do not have torsion are the trivial knot, the Hopf link, their connected sums and disjoint unions. The nontrivial torsion always contains the $Z_2$ subgroup.
  \vskip .3cm
  
  They further conjecture (and prove in many cases) that the above theorem will also hold for any, not necessarily alternating links.
  
  The above statements allow us to reformulate the geometric definition in terms of torsion.  But first we have to eliminate the case of the Hopf link.
    \vskip .3cm
  
   {\bf Calculation}. The Khovanov homology of the Hopf link is quasiisomorphic to $(Z_2, 0, Z_2)$.
 
    \vskip .3cm
    
    Assigning 0 and 1-resolutions to 2 crossing points of the Hopf link one gets 4 complete resolutions
    and the Khovanov complex becomes:
    
    $$0 \rightarrow V  \otimes V \rightarrow V \oplus V \rightarrow V \otimes V \rightarrow 0$$
    
    recall that that the space $V$ is generated by $v_+$ and $v_-$. One can see that the 1-cycles are
    $(v_- \otimes v_-)$, $(v_+ \otimes v_- - v_- \otimes v_+)$ and 1-boundaries are 0. The 2-cycles are 
    $(v_+^1 \oplus v_+^2)$ and $(v_-^1 \oplus v_-^2)$, 2-boundaries are $(v_+^1 \oplus v_+^2)$ and 
    $(v_-^1 \oplus v_-^2)$, 3-cycles are $(v_+ \otimes v_+)$, $(v_+ \otimes v_- )$ , $ (v_- \otimes v_+)$,
    3-boundaries are $(v_- \otimes v_-)$, $(v_+ \otimes v_-  + v_- \otimes v_+)$. Thus the Khovanov
    homology of the Hopf link is  $(Z_2, 0, Z_2)$.
        \vskip .3cm
    
  Consider a new invariant,  an additive functor:
  
  $$\mathcal T_2: C \in Ob(\mathcal D) \rightarrow {Tor}_2( H^*(C))$$
    \vskip .3cm
  {\bf Definition (T)}. The local system is of finite type n if for any codimension $n$ selfintersection of the discriminant the corresponding nth cone is not quasiisomorphic to $(Z_2, 0, Z_2)$ and has  torsion-free homology, i.e. the image of  $\mathcal T_2$
  is zero.
    \vskip .3cm  
 The disadvantage of this definition is that after taking torsion the theory becomes trivial on the level
 of the Euler characteristics.  To fix this problem, we consider the ideals $ I_n$ generated by the complexes $C(n)^{\bullet}$  introduced in section 6.

  \vskip 1.5cm
 \centerline {\bf 8. The Filtration.}

\vskip 1cm

 In this paragraph we take the factorization, which we discussed in section 5, to the level of the derived category.

According to [D], it is possible to divide in the tensor monoidal categories if the tensor product preserves quasiisomorphisms, i.e. an exact functor.  If we consider the class of flat objects (those, on which the tensor product is exact), then one can divide by a subcategory.

 Recall our geometric definition of finiteness:
 \vskip .3cm   
{\bf Definition (G)}. The local system is of finite type n, if there exists such minimal $n$, s.t. for any selfintersection of the discriminant of codimension $n$ the corresponding complex is quasiisomorphic to
$C^{\bullet}(n)(\mathcal U) $, where $\mathcal U$ is the Khovanov complex of the disjoint union of unknots.
  \vskip .3cm

   We will follow the Verdier approach, who constructed the quotient 
   $\mathcal T/ \mathcal S$, where $ \mathcal S$ is a subcategory of $\mathcal T$. He showed that the factorization is well-defined if $ \mathcal S$ is thick.
    Recall the definition of a {\bf thick} subcategory:
   \vskip .3cm  
  
  {\bf Definition}. The subcategory of a triangulated category is called {\bf thick} if it is triangulated and contains all direct summands of it's objects.  
 \vskip .3cm
   
   The subcategories, over which we will be factorizing are generated by $ I_n$ and {\bf supported on $D_n$} - the union of strata of codimension $n$ of the discriminant:
      \vskip .3cm

First we define the notion of generated subcategory:

 \vskip .3cm
 {\bf Definition}.  Let $\mathcal T$ be a triangulated category satisfying the {\bf homotopy colimits} axiom.
 Let $\alpha$ be an infinite cardinal. Let $S$ be a class of objects of $\mathcal T$. Then 
 ${<S>}^{\alpha}$
 will denote the smallest  $ \mathcal S$,  $ \mathcal S$ a triangulated subcategory of  $ \mathcal T$, the {\bf generated subcategory} satisfying:
 
 1). The objects of S lie in  $ \mathcal S$.
 
 2). Any coproduct of fewer that $\alpha$ objects of  $ \mathcal S$ lies in  $ \mathcal S$.
 
 3). The subcategory  $ \mathcal S \subset   \mathcal T$ is thick.
  \vskip .3cm

 We refer to [N] for the proofs that  $ \mathcal S$ is well-defined and that it is localizing [Bo]. By Verdier theorem [V] one can factorize by such subcategories.
   \vskip .3cm
      
    {\bf Proposition 5}. Let ${\mathcal I_n}$ be the  triangulated category generated by 
    $C^{\bullet}(n)$. Then for any complex $X^{\bullet} \in {\mathcal I_n}$ we have $\chi (X^{\bullet})=2^n \cdot k$.
      \vskip.2cm 
    {\bf Proof}. This is a check for all operations in the triangulated category:
    
    1) quasiisomorphic complexes have the same Euler characteristics
    
    2) By taking the direct sum of  i copies of $C^{\bullet}(n)$ we get a complex with Euler characteristics
    $2^n \cdot i$ :
    
     $$\chi (\underbrace{X^{\bullet} \oplus X^{\bullet} \oplus... X^{\bullet}}_{i times}) = i \chi(X^{\bullet})$$
    
    3) By taking the  cone of the map between complexes $C^{\bullet}(n)$ we can arrange any combination of morphisms between the components, but it is an easy check that what we get will 
 have    a trivial Euler characteristics. 
    
    In all cases the Euler characteristics of the resulting complex will be a multiple of $2^n$. Thus the divisability of the Euler characteristics by $2^n$ becomes an invariant of $\mathcal I_n$

     \vskip.2cm 
   {\bf Proposition 6.}  One gets a sequence of  categories:
        $$ ...{\mathcal I_n} \subset  {\mathcal I_{n-1}} \subset ... \subset  {\mathcal I_1} $$ 
   
        \vskip.2cm 
       {\bf Proof}. This statement follows from the previous one: if one takes k  equal 2, and recall that
       ${\mathcal I_n}$ is supported on $D_n$, then
       $${\mathcal I_n} \subset  {\mathcal I_{n-1}} $$

    {\bf Definition }.      We define the derived category $\tilde{\mathcal D_n}$ as a factor category (in a sense of Verdier [V], [D]) : 
  
  $$\tilde{\mathcal D_n} = \tilde{\mathcal D}/ \mathcal I_n$$
  
  \vskip .3cm
   Where  $I_n$ is a thick subcategory.

 It follows from the above definition that the subcategory is thick when it is closed under cofibrations, retractions, direct sums and
suspensions.
 When we form a thick subcategory, generated by  $\mathcal I_n$, we may loose the property of the Euler characteristics being divisible by $2^n$, but since $\mathcal I_n$ is {\bf supported} on $D_n$,
 this is the thick subcategory of $\tilde{\mathcal D}$. The resulting category $\tilde{\mathcal D_n}$ will
 lack the exactness property only on codimension $n$ strata.

   \vskip .3cm
   
   {\bf Example.} In the triangulated category of complexes over $Z$ consider
   those, which have $k$-torsion in the homology.
   This subcategory is thick.
    \vskip .3cm
    {\bf Proposition 7}. The Verdier quotients of the category $\tilde{\mathcal D}$ over thick subcategories 
    $\tilde{\mathcal{I}_n}$ defined above, form an increasing filtration of the category $\tilde{\mathcal D}$:
    
      $$ \tilde{\mathcal D}=\tilde{\mathcal D_{\infty}} \supset ...\tilde{\mathcal D_n} \supset  \tilde{\mathcal D_{n-1}} \supset ... \supset  \tilde{\mathcal D_1} $$
    
   \vskip.2cm 
       {\bf Proof}.  The statement follows from the previous propositions.

   \vskip .5cm

  {\bf Remark 1.} Notice that in  $ \tilde{\mathcal D}$ the Euler characteristics of complexes is well-defined
  modulo  $2^n$.
  
  \vskip .5cm

 {\bf Remark 2}. Notice that this final definition is very similar to the one of invariants of finite type given in paragraph 3.
 \vskip .5cm

   {\bf Remark 3}. Yet another approach...  Complexes with torsion-free homology don't form a thick subcategory, but they are factor- functors of Tor. We show that one still can form a factor-category $\tilde{\mathcal D_{tor}}$ with the same objects, morphisms (arrows) of which can factor through complexes, homology of which don't have torsion [Ho].
    
     \vskip .3cm
     
  {\bf Theorem}. In the factor-category $\tilde{\mathcal D_{tor}}$ objects, corresponding to complexes, homology of which have no torsion are isomorphic to zero, while complexes, homology of which have torsion are never isomorphic to zero in $\tilde{\mathcal D_{tor}}$.

      \vskip.2cm 
   {\bf Definition (V')}. The local system of Khovanov complexes is of finite type $n$ if for any codimension $n$ selfintersection of the discriminant the nth cone is zero in $\tilde{\mathcal D_{tor}}$.

\vskip 1.5cm

 \centerline {\bf 9. The Finiteness result.}

\vskip 1cm 
  In this paragraph we prove the first simple finiteness property
  of  the Khovanov local system. 
  
   First let's define our categories.
   
   The category of knots $\mathcal K$ is the topolocical category, objects of which are knots and
   morphism are knot cobordisms. ( If knots $K_1$ and $K_2$ have the same isotopy type, i.e.
   lie within the same chamber of the Vassiliev space, then morphisms are just the product cobordisms,
   if we pass between adjacent chambers, changing one crossing, then these are genus one cobordisms). 
   
   The construction of the local system of Khovanov complexes on the Vassiliev space of knots [S1] provided us with a functor from
   the category of knots into the derived category of complexes, denote it by $\mathcal Kh$.

    Recall that the {\bf crossing number} of the knot is the minimum of the crossing numbers over all
    its projections.
    
     In the category of knots $\mathcal K$  we  define a sequence of subcategories $K_n$ . Objects
    of $\mathcal K_n$ are knots with at most n crossings, morphisms in $\mathcal K_n$  are 
    cobordisms between knots with the crossing number at most n, etc. The corresponding derived
    category is denoted $\mathcal Kh_n$.

 \vskip .3cm
 {\bf Theorem 1}. Restricted to the subcategory of knots with at most $n$ crossings , $n\geq 3$, Khovanov local system is of finite type $\leq n$.

   \vskip .5cm

{\bf Proof}. 
We give the proof according to the  geometric definitions.

\vskip .3cm

If we restrict Khovanov theory to the subcategory $\mathcal K_n$ of knots  with at most $n$ crossing, 
the complexes we get from knot projections are all quasiisomorphic to the ones of length at most $(n+1)$.

Consider an n-dimensional commutative hypercube, corresponding to the selfintersection of the discriminant
of codimension n. 

As we have seen in the previous paragraph, complexes, corresponding to the codimension one  walls of the discriminant are

$$
\renewcommand\arraystretch{1}
X^{\bullet}\oplus X^{\bullet}[2] \quad\mathrm{with \ the \ matrix \ differential}\quad
d_{C_\omega}=\left(\!\!\!\begin{array}{cc}
d_X & 0 \\ 0 & d_X[2] \end{array}\!\!\right),
$$
where $X^{\bullet}$ is the Khovanov complex of the knot of (n-1) crossings where  kth double point of $K$ is given 1-resolution.

 When we pass to codimension 2 selfintersections of the discriminant, we get 4-graded complexes associated to it:

$$Y^{\bullet} \oplus Y^{\bullet}[2] \oplus Y^{\bullet}[2] \oplus Y^{\bullet}[4]$$

\vskip .3cm
where $Y^{\bullet}$  is the Khovanov complex of the knot of (n-2) crossings where  kth and lth double points of $K$ are given 1-resolution.

After establishing  these identities for all indices up to $n$, it is easy to see from our geometrical interpretation
of the finite type condition what is the nth generalized cone of the restricted Khovanov local
system.

 The geometric definition implies that  we calculate the nth cone by taking the "first" knot in the hypercube (via the coorientation)	and giving all the crossings of it's projection
1-resolutions. We end up with a collection of circles in $R^2$, corresponding to the last component of
the complex, associated with the projection of the knot.

 The $n$th cone, assigned to the selfintersection of the discriminant of codimension $n$ will be
 quasiisomorphic to the complex of the disconnected sum of $2^n$ copies of  the collection of circles, described above, shifted in homological grading according to the coorientation.	 In the stable category this complex is isomorphic to the one of the circle.

 Passing to the stable category means the eliminating of the "inessential" crossings of the knot projection.  The discriminant of the space of knots, consists to singular knots with "essential" crossings (i.e. this singular knot can be realized without extra under/over crossings). Let $\mathcal K_n$ be the subcategory of knots with the crossings number (minimum over all projections) at most $n$ , then we can reformulate the main theorem as follows:

\vskip .3cm
{\bf Theorem 1'}.  If the singular knot $K$ with n double points can be realized without extra under/over crossings, then the extension of the Khovanov local system to $K$ is zero in $\tilde{\mathcal D_n}$  (even if in the diagram of
$K$  there are other under/over crossings).
\vskip.3cm
Next we consider the localized local system: since the nth cone is given by the formula $C^\bullet (n) \otimes K_{1,1,....1}$, where $K_{1,1,....1}$ Khovanov complex of the knot projection, where n crossing points are given 1-resolutions and

\[ C^\bullet (n) = 
\left (  
	\begin{array}{c}
   		{C_n^n}  Z \\ 
		0 \\ 
		\vdots\\ 
		0 \\ 
		{C^{i}_{n}}  Z \\ 
		0 \\ 
		\vdots\\ 	
		0 \\
		{C_n^0}  Z  
	\end{array} 
\right) \]

in  $\tilde{\mathcal D_n}$ it will satisfy the finiteness condition, since all crossing points will be given 1- resolution and this is a local system of finite type n.
And since we are not considering walls with knots of crossing number $n+1$, this is the only cone we can form.
 \vskip .5cm
 
 {\bf Remark}. One may not need to take the nth cone to get the finiteness condition, i.e. the knot projection may decompose into a disconnected sum of unknots earlier. (It would be interesting to get an estimate). But since $\tilde{\mathcal D_n} \supset  \tilde{\mathcal D_k}$ for $k \leq n$, we get that the local system will be of type at most $n$.
 \vskip .5cm

  {\bf Example}.  Consider the right-handed trefoil with three double points:

   There are 8 resolutions of this singular knot, corresponding to 8 chambers (say, complexes X,Y,A,B,C,D,W,Z), adjacent to the selfintersection of the discriminant of codimansion 3. Two of them correspond to the righthanded trefoil and it's mirror image and six - to the twisted unknots, which are obtained after changing any overcrossing in the projection of the trefiol  to the undercrossing.   
  \vskip 1cm

\begin{equation*}
\xymatrix
{& & & & 0 \ar[d] &   \\
& & 0  \ar[d]^{d_Y} & \ar[r]^{\omega}  & { Z_0 = X_3} \ar[d]^ {d_Z}   &  \\
& 0 \ar[d]^{d_X} & {Y_0} \ar[d]^{d_Y}  \ar[r]^{\omega} & {...}  \ar[r]^ {\omega}  & {...}  \ar[d]^{d_Z} &  \\
& {X_0} \ar[d]^{d_X}  \ar[r]^{\omega} & {Y_1} \ar[d] ^{d_Y}  & {...} & {Z_3 = X_0}  \ar[d] & \\
& {X_1} \ar[d]^{d_X}  \ar[r]^{\omega} & {Y_2}  \ar[d] ^{d_Y}    \ar[r]^ {\omega} & {...} & 0 & \\
 & {X_2}  \ar[d]^{d_X}  \ar[r]^{\omega} & {Y_3} \ar[d] ^{d_Y}  \ar[r]^{\omega} & & & \\
 &{X_3} \ar[d] \ar[r]^{\omega}  & 0 & & & \\
& 0 & & & & \\ }
\end{equation*}

  \vskip 1cm 
  
  The third cone $C^\bullet (3) $ will be of the form:
  
 \[ C^\bullet (3) \otimes K_{1,1,1} =  \left (  \begin{array}{c}
   Z \\ 
   0 \\
    3 Z \\
     0 \\
      3 Z \\
       0 \\ 
        Z \end{array} \right )  \otimes K_{1,1,1} \]

 \vskip .5cm 
 
 Thus the third cone is acyclic in $\tilde{\mathcal D_3}$ and this is a condition for the local system to be of type 3.

 In the upcoming paper [S2] we prove the generalization of Theorem 1, or the categorification of Birman-Lin theorem.

\vskip 1.5cm 

\newpage
 \centerline {\bf 10. Further directions.}
\vskip 1cm 

1. We would like to generalize  the  result  of this paper to the Khovanov-Rozhansky homology [KR],
for which we defined the wall-crossing morphisms [SW].

\vskip .3cm
2. It would be very interesting to see, if the recent extensions to the singular locus of Ozsvath-Szabo knot invariants, done
by Benjamin Audoux [A] , cf.  [OSS],  satisfy finiteness conditions.

\vskip .3cm

3. The geometric definition implies that the knot homology theory is of finite type, if after taking sufficiently high cones all objects decompose into a collection of disconnected circles. It would be very interesting to understand what are the "building blocks" for the homology theories of higher dimensions. 
\vskip .3cm
4. Our theorem can be proved using the properties of the homological width of the knot. Such proof
could provide another point of view on finiteness result (via the ranks of homologies).
\vskip .3cm
5. We believe that our constructions will provide a better estimates on the crossing number of the knot.

\vskip 3cm

nadya@math.stanford.edu

\newpage
 \centerline {\bf 11. Bibliography.}
\vskip 1cm 
[A] Adams J. F., Stable homotopy and generalized homology, Chicago Lectures in Mathematics, Univ. of Chicago Press, 1974.
\vskip .2cm
[AP] Asaeda M., Przytycki J., Khovanov homology: torsion and thickness., arXiv:0402402.
\vskip .2cm
[Au] Audoux B., Heegaard-Floer homology for singular knots, arcXiv:0705.2377.
\vskip .2cm

[B]  Bar-Natan D.,Khovanov's homology for tangles and cobordisms,  arXiv:math/0410495.
\vskip .2cm
[Bo] Bousfield A. The localization of spaces with respect to homology, Topology 14, 1975, 133-150.
\vskip .2cm
[BL]  Birman J., Lin X.S., Knot polynomials and Vassiliev invariants, Invent. Math., 111 (1993), pp.225-270.
\vskip .2cm
[D] Drinfeld D., DG quotients of DG categories, arXiv:0210114v6
\vskip .2cm
[GM] Gelfand S., Manin Yu., Methods of homological algebra, Springer 1996.
\vskip .2cm
[H] Hovey M., Model categories, Math.Surveys, vol.63.
\vskip .2cm
[Kh] Khovanov M., A Categorification of the Jones Polynomial,  Duke Math. J. 101 (2000), no. 3, 359--426.

\vskip .2cm
[Kh1] Khovanov M., Patterns in knot cohomology,  I, Experiment. Math. 12 (2003), no. 3, 365-374 math.QA/0201306.
\vskip .2cm
[Kh2] Khovanov M., A functor-valued invariant of tangles., arXiv:math/0103190.
\vskip .2cm

[Kh3] Crossingless matchings and the cohomology of $(n,n)$ Springer varieties, arXiv:QA/0202110. 
\vskip .2cm
[KR] Khovanov M., Rozansky L., Matrix factorizations and link homology
\vskip .2cm
[N] Neeman A., Triangulated categories, Princeton university press, 2001.
\vskip .2cm
[OSS]  Ozsvath P.,  Stipsicz A., Szabo Z., Floer homology and singular knots,  arXiv:0705.2661.
\vskip .2cm

\vskip .2cm
[P]  Przytycki J., When the theories meet: Khovanov homology as Hochschild homology of links, 
arXiv GT/0509334.
\vskip .2cm
[Q] Quillen D., Homotopical Algebra. Lecture Notes in Math., no. 43, Springer-Verlag, 1967
\vskip .2cm
[S1] Shirokova N., On the classification of Floer-type theories, arXiv:0704.1330.
\vskip .2cm
[S2] Shirokova N., The Categorification of Birman-Lin Theorem, preprint 2007.

\vskip .2cm
[SW] Shirokova N., Webster B., Wall-crossing morphism for Khovanov-Rozhansky homology, arXiv:0706.1388.
\vskip .2cm

[Su] Shumakovich A., Torsion of the Khovanov homology, arXiv:GT/0405474v1.
\vskip .2cm
[V] Verdier J. L., Des categories derivees des categories abeliennes. Asterisque, vol.239. Societe mathematique de France, 1996.
\vskip .2cm

[Vo] Voevodsky V., The $A^1$-homotopy theory, Documenta Math. J.DMV, 1998, I,579-604.
\vskip .2cm

\end{document}